\documentclass[conference]{IEEEtran}

\usepackage{cite}
\usepackage[dvips]{graphicx}
\usepackage[cmex10]{amsmath}
\usepackage{algorithmic}
\usepackage{algorithm}
\usepackage{amssymb}
\usepackage{balance}
\usepackage{subfigure}
\usepackage{cases}

\DeclareMathOperator*{\argmin}{arg\,min}
\DeclareMathOperator*{\x}{\mathbf{x}}
\DeclareMathOperator*{\y}{\mathbf{y}}
\DeclareMathOperator*{\z}{\mathbf{z}}
\DeclareMathOperator*{\A}{\mathbf{A}}

\DeclareMathOperator*{\blambda}{\boldsymbol \lambda}
\DeclareMathOperator*{\R}{\mathbb{R}}

\begin{document}
\title{Multi-Block ADMM for Big Data Optimization\\ in Modern Communication Networks}

\author{\IEEEauthorblockN{Lanchao Liu and Zhu Han}
\IEEEauthorblockA{Department of Electrical and Computer Engineering\\
University of Houston,
Houston, TX, 77004}}


\maketitle

\begin{abstract}
In this paper, we review the parallel and distributed optimization algorithms based on the alternating direction method of multipliers (ADMM) for solving ``big data'' optimization problems in modern communication networks. We first introduce the canonical formulation of the large-scale optimization problem. Next, we describe the general form of ADMM and then focus on several direct extensions and sophisticated modifications of ADMM from $2$-block to $N$-block settings to deal with the optimization problem. The iterative schemes and convergence properties of each extension/modification are given, and the implementation on large-scale computing facilities is also illustrated. Finally, we numerate several applications in communication networks, such as the security constrained optimal power flow problem in smart grid networks and mobile data offloading problem in software defined networks (SDNs).
\end{abstract}

\section{Introduction}
Nowadays, modern communication networks play an important role in electric power system, mobile cloud computing, smart city evolution and personal health care. The employed novel telecommunication technologies make data collection much easier for power system operation and control, enable more efficient data transmission for mobile applications, and promise a more intelligent sensing and monitoring for metropolitan city-regions. Meanwhile, we are witnessing an unprecedented rise in volume, variety and velocity of information in modern communication networks. A large volume of data are generated by our digital equipments such as mobile devices and computers, smart meters and household appliances, as well as surveillance cameras and sensor-equipped mass rapid transit around the city. The information exposition of big data in modern communication networks makes statistical and computational methods significantly important for data analysis, processing, and optimization. The network operators or service providers who can develop and exploit efficient methods to tackle big data challenges will ensure network security and resiliency, gain market share, increase revenue with distinctive quality of service, as well as achieve intelligent network operation and management.

The unprecedented ``big data'', reinforced by communication and information technologies, presents us opportunities and challenges. On one hand, the inferential power of algorithms, which have been shown to be successful on modest-sized data sets, may be amplified by the massive dataset. Those data analytic methods for the unprecedented volumes of data promises to personalized business model design, intelligent social network analysis, smart city development, efficient healthcare and medical data management, and the smart grid evolution. On the other hand, the sheer volume of data makes it unpractical to collect, store and process the dataset in a centralized fashion. Moreover, the massive datasets are noisy, incomplete, heterogeneous, structured, prone to outliers, and vulnerable to cyber-attacks. The error rates, which are part and parcel of any inferential algorithm, may also be amplified by the massive data. Finally, the ``big data" problems often come with time constraints, where a medium-quality answer that is obtained quickly can be more useful than a high-quality answer that is obtained slowly. Overall, we are facing a problem in which the classic resources of computation such as time, space, and energy, are intertwined in complex ways with the massive data resources.

With the era of ``big data'' comes the need of parallel and distributed algorithms for the large-scale inference and optimization. Numerous problems in statistical and machine learning, compressed sensing, social network analysis, and computational biology formulates optimization problems with millions or billions of variables. Since classical optimization algorithms are not designed to scale to problems of this size, novel optimization algorithms are emerging to deal with problems in the ``big data'' setting. An incomprehensive list of such kind of algorithms includes block coordinate descent method \cite{PS09,YS09,Y12}\footnote{\cite{Y12} proposes a stochastic block coordinate descent method.}, stochastic gradient descent method \cite{LO08, MMAL10, FBCS11}, dual coordinate ascent method \cite{CKCSS08, ST13}, alternating direction method of multipliers (ADMM) \cite{BT97, BPCPE10} and Frank-Wolf method (also known as the conditional gradient method) \cite{RP14,SMMP13}. Each type of the algorithm on the list has its own strength and weakness. The list is sill growing and due to our limited knowledge and the fast develop nature of this active field of research, many efficient algorithms are not mentioned here.

In this paper, we focus on the application of ADMM for the ``big data'' optimization problem in communication networks like smart grids and software defined networks (SDNs). In particular, we consider the parallel and distributed optimization algorithms based on ADMM for the following convex optimization problem with a canonical form as
\begin{align}
\label{eqn:intro}
\min_{\x_1,\x_2, \ldots, \x_N}\quad & f(\x) = f_i({\x}_i) + \ldots + f_i({\x}_N),\nonumber \\
\text{s.t.}\quad &{\A}_i{\x}_i + \ldots + {\A}_N{\x}_N  = \mathbf{c}, \nonumber \\
&{\x}_i \in \mathcal{X}_i, \quad i = 1,\ldots,N,
\end{align}
where $\x = ({\x}_1^{\top}, \ldots, {\x}_N^{\top})^{\top}$, $\mathcal{X}_i \subset {\R}^{n_i} (i = 1,2, \ldots, N)$ are closed convex set, ${\A}_i \in {\R}^{m \times n_i}(i = 1,2, \ldots, N)$ are given matrices, $\mathbf{c} \in {\R}^{m}$ is a given vector, and $f_i: {\R}^{n_i} \rightarrow {\R}$ $(i = 1,2, \ldots, N)$ are closed convex proper but not necessarily smooth functions, where the non-smoothness functions are usually employed to enforce structure in the solution. Problem (\ref{eqn:intro}) can be extended to handle linear inequalities by introducing slack variables. Problem (\ref{eqn:intro}) finds wide applications in smart grid on distributed robust state estimation, network energy management and security constrained optimal power flow problem, which we will illustrated later.

Though many algorithms can be applied to deal with problem (\ref{eqn:intro}), we restrict our attention to the class of algorithms based on ADMM. The rest of this paper is organized as follows. Section \ref{sec:background} introduces the background of the ADMM and its two direct extensions for problem (\ref{eqn:intro}) to $N$ blocks. The limitations of those direct extensions are also addressed. Section \ref{sec:multiblock} gives three approaches based on Variable Splitting, ADMM with Gaussian back substitution and proximal Jacobian ADMM to the multi-block settings, respectively, for problem (\ref{eqn:intro}) with provable convergence. The applications of problem (\ref{eqn:intro}) in communication networks are described in Section \ref{sec:application}. Specifically, we discuss two examples in detail: the security constrained optimal power flow problem in smart grid networks and mobile data offloading problem in SDNs. Section \ref{sec:conclusion} summarizes this paper.

\section{ADMM Background}
\label{sec:background}

In this section, we first introduce the general form of ADMM for optimization problem analogous to (\ref{eqn:intro}) with only two blocks of functions and variables. After that, we describe two direct extensions of ADMM to multi-block setting.

\subsection{ADMM}
The ADMM was proposed in \cite{RA75}, \cite{DB76} and recently revisited by \cite{BPCPE10}. The general form of ADMM is expressed as
\begin{equation}
\label{eqn:StdAMDD}
\min_{\x_1 \in \mathcal{X}_1,\x_2 \in \mathcal{X}_2} f_1({\x}_1) + f_2({\x}_2) \quad \text{s.t.} \quad {\A}_1{\x}_1 + {\A}_2{\x}_2 = \mathbf{c}.
\end{equation}
The augmented Lagrangian for (\ref{eqn:StdAMDD}) is
\begin{align}
\mathcal{L}_{\rho} ({\x}_1,{\x}_2,{\blambda})  &= f_1({\x}_1) + f_2({\x}_2) - {\blambda}^{\top}( {\A}_1{\x}_1 + {\A}_2{\x}_2 - \mathbf{c}) \nonumber \\
&+ \frac{\rho}{2}\Vert {\A}_1{\x}_1 + {\A}_2{\x}_2 - \mathbf{c} \Vert_2^2,
\end{align}
where $\blambda \in {\R}^{m}$ is the Lagrangian multiplier and $\rho > 0$ is the parameter for the quadratic penalty of the constraints. The iterative scheme of ADMM embeds a Gauss-Seidel decomposition into iterations of ${\x}_1$ and ${\x}_2$ as follows

\begin{numcases}{}
\label{eqn:update}
{\x}_1^{k+1} = \argmin_{{\x}_1}\mathcal{L}_{\rho} ({\x}_1,{\x}_2^{k},{\blambda}^{k}),\\
\label{eqn:update_2}
{\x}_2^{k+1} = \argmin_{{\x}_2}\mathcal{L}_{\rho} ({\x}_1^{k+1},{\x}_2,{\blambda}^{k}),\\
{\blambda}^{k+1} = {\blambda}^{k} - {\rho}({\A}_1{\x}_1^{k+1} + {\A}_2{\x}_2^{k+1} - \mathbf{c}),
\end{numcases}
where in each iteration, the augmented Lagrangian is minimized over ${\x}_1$ and ${\x}_2$ separately. In (\ref{eqn:update}) and (\ref{eqn:update_2}), functions $f_1$ and $f_2$ as well as variables ${\x}_1$ and ${\x}_2$ are treated individually, so easier subproblems can be generated. This feature is quite attractive and advantageous for a broad spectrum of applications. The convergence of ADMM for convex optimization problems with two blocks of variables and functions has been proved in \cite{BT97}, \cite{BPCPE10},  and the iterative scheme is illustrated in Algorithm \ref{alg:A1}. Algorithm \ref{alg:A1} can deal with multi-block case when auxiliary variables are introduced, which will be described in Section \ref{sec:VSADMM}.

\begin{algorithm}[t]
\caption{Two-block ADMM}\label{alg:A1}
\begin{algorithmic}
\STATE Initialize: ${\x}^{0}$, ${\blambda}^{0}$, $\rho>0$;

\FOR{$k=0,1,\ldots$}

\STATE ${\x}_1^{k+1} = \argmin_{{\x}_1}\mathcal{L}_{\rho} ({\x}_1,{\x}_2^{k},{\blambda}^{k})$;

\STATE ${\x}_2^{k+1} = \argmin_{{\x}_2}\mathcal{L}_{\rho} ({\x}_1^{k+1},{\x}_2,{\blambda}^{k})$;

\STATE ${\blambda}^{k+1} = {\blambda}^{k} - {\rho}({\A}_1{\x}_1^{k+1} + {\A}_2{\x}_2^{k+1} - \mathbf{c})$;

\ENDFOR

\end{algorithmic}
\end{algorithm}

\subsection{Direct Extensions to Multi-block Setting}
\begin{figure}[t]
    \centering
    \includegraphics[width=0.48\textwidth]{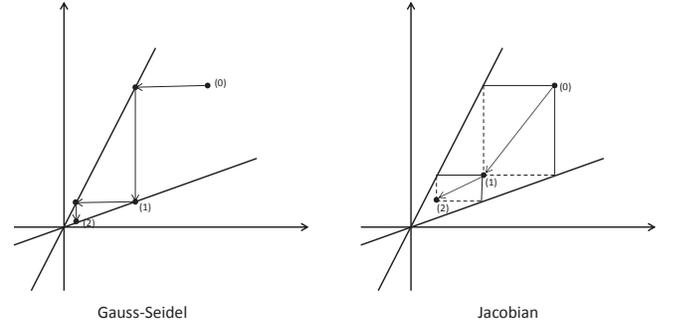}
    \caption{Comparison of Gauss-Seidel update and Jacobian update.}
    \label{fig:GSJ}
\end{figure}

The ADMM promises to solve the optimization problem (\ref{eqn:intro}) with the same philosophy as Algorithm \ref{alg:A1}. In the following, we present two kinds of direct extensions, Gauss-Seidel and Jacobian for multi-block ADMM. The comparison of these two updates is shown in Figure \ref{fig:GSJ}, To be specific, we first give the augmented Lagrangian function of problem (\ref{eqn:intro})
\begin{align}
\label{eqn:ALM}
\mathcal{L}_{\rho} ({\x}_1,\ldots,{\x}_N,{\blambda})  = &\sum_{i=1}^{N}f_i({\x}_i) - {\blambda}^{\top}( \sum_{i=1}^{N} {\A}_i{\x}_i - \mathbf{c})\\ \nonumber
&+ \frac{\rho}{2}\Vert \sum_{i=1}^{N} {\A}_i{\x}_i - \mathbf{c} \Vert_2^2.
\end{align}
where a penalty on linear constrains is added to the Lagrangian function. The $\rho$ is the penalty parameter.

\subsubsection{Gauss-Seidel}
Intuitively, a natural extension of the classical Gauss-Seidel setting ADMM from $2$ blocks to $N$ blocks is a straightforward replacement of the two-block alternating minimization scheme by a sweep of update of ${\x}_i$ for $i = 1,2,\ldots,N$ sequentially. In particular, at iteration $k$, the update scheme for ${\x}_i$ is
\begin{equation}
{\x}_i = \argmin_{{\x}_i}\mathcal{L}_{\rho}(\{{\x}_j^{k+1}\}_{j<i},{\x}_i,\{{\x}_j^{k}\}_{j>i},{\blambda}^k),
\end{equation}
where $\{{\x}_j\}_{j<i}$ denotes the set of variables prior to $i$. The augmented Lagrangian function (\ref{eqn:intro}) is split and updated alternatingly. The direct Gauss-Seidel type extension can be illustrated in Algorithm \ref{alg:A2}.

\textbf{Remark:} Algorithm \ref{alg:A2} has been utilized in practical problems \cite{YAJWY12, MX11, HCB14} despite a lack of rigourous proof for the convergence. Actually, the convergence of Gauss-Seidel multi-block ADMM is not well understood and is ambiguous for a long time: Neither affirmative convergence proof nor counter examples for convergence failure are shown in the literature. Recently, \cite{CBYX13} has shown that the direct extension of Gauss-Seidel mulit-block ADMM is not necessarily convergent. \cite{MZ12} proves the convergence of Algorithm \ref{alg:A2} with a sufficient small step size for Lagrangian multiplier update and additional assumptions on the problem (\ref{eqn:intro}). \cite{CBYX14} conjectures that an independent uniform random permutation of the update order for blocks in each iteration will result in a convergent iteration scheme.  \cite{BXM12,MTXMSZ14} propose some slightly modified version of Algorithm \ref{alg:A2} with provable convergence and competitive iteration simplicity and computing efficiently, which we will illustrate later in Section \ref{sec:GSADMM}.

\begin{algorithm}[t]
\caption{Gauss-Seidel Multi-block ADMM}\label{alg:A2}
\begin{algorithmic}
\STATE Initialize: ${\x}^{0}$, ${\blambda}^{0}$, $\rho>0$;

\FOR{$k=0,1,\ldots$}

\FOR{$i=1,\ldots,N$}

\STATE \COMMENT{${\x}_i$ is updated \textbf{sequentially}.}

\STATE ${\x}_i^{k+1} = \argmin_{{\x}_i}\mathcal{L}_{\rho}(\{{\x}_j^{k+1}\}_{j<i},{\x}_i,\{{\x}_j^{k}\}_{j>i},{\blambda}^k)$;

\ENDFOR

\STATE ${\blambda}^{k+1} = {\blambda}^{k} - {\rho}(\sum_{i=1}^{N} {\A}_i{\x}_i^{k+1} - \mathbf{c})$;

\ENDFOR

\end{algorithmic}
\end{algorithm}

\subsubsection{Jacobian}
Another possible iterative scheme for the $N$ blocks ADMM is the Jacobian type update, which performs the update of ${\x}_i$ in a parallel coordinate fashion for $i = 1, \ldots, N$. In particular, the update of ${\x}_i$ is calculated as:
\begin{equation}
{\x}_i = \argmin_{{\x}_i}\mathcal{L}_{\rho}({\x}_i,\{{\x}_j^{k}\}_{j \neq i},{\blambda}^k),
\end{equation}
where $\{{\x}_j^{k}\}_{j \neq i}$ denotes the set of variables except for ${\x}_i$. Different from the iterative scheme of Algorithm \ref{alg:A2} that the update of ${\x}_i$ has to be performed sequentially one after another, the iterations in the Jacobian ADMM can be performed concurrently, i.e. all ${\x}_i$ can be updated in a parallel fashion. This advantage makes the Jacobian type ADMM preferred for parallel implementation, and the direct Jacobian type extension can be illustrated in Algorithm \ref{alg:A3}.

\textbf{Remark} Though Algorithm \ref{alg:A3} is more computational efficient in the sense of parallelization, \cite{BLX13} shows that Algorithm \ref{alg:A3} is not necessarily convergent in the general case, even in the 2 blocks case. \cite{WMZW14} proves that if matrices ${\A}_i$ are mutually near-orthogonal and have full column-rank, Algorithm \ref{alg:A3} converges globally. A proximal Jacobian ADMM is also proposed in \cite{WMZW14} with provable convergence, which we will illustrate later in Section \ref{sec:PJADMM}.

\begin{algorithm}[t]
\caption{Jacobian Multi-block ADMM}\label{alg:A3}
\begin{algorithmic}
\STATE Initialize: ${\x}^{0}$, ${\blambda}^{0}$, $\rho>0$;

\FOR{$k=0,1,\ldots$}

\FOR{$i=1,\ldots,N$}

\STATE \COMMENT{${\x}_i$ is updated \textbf{concurrently}.}

\STATE ${\x}_i^{k+1} = \argmin_{{\x}_i}\mathcal{L}_{\rho}({\x}_i,\{{\x}_j^{k}\}_{j \neq i},{\blambda}^k)$;

\ENDFOR

\STATE ${\blambda}^{k+1} = {\blambda}^{k} - {\rho}(\sum_{i=1}^{N} {\A}_i{\x}_i^{k+1} - \mathbf{c})$;

\ENDFOR

\end{algorithmic}
\end{algorithm}

\section{Multi-block ADMM}
\label{sec:multiblock}
In this section, we introduce several sophisticated modifications of ADMM, Variable splitting ADMM \cite{BT97, BPCPE10, MJM10}, ADMM with Gaussian Back Substitution \cite{BXM12,BMX12} and Proximal Jacobian ADMM \cite{WMZW14,BMX13}, to deal with the multi-block setting.
\subsection{Variable Splitting ADMM}
\label{sec:VSADMM}
To solve the optimization problem (\ref{eqn:intro}), we can apply the variable splitting \cite{BT97, BPCPE10, MJM10} to deal with the multi-block variables. In particular, the optimization problem (\ref{eqn:intro}) can be reformulated by introducing auxiliary variable $\z$
\begin{align}
\label{eqn:VS}
\min_{\x,\z}\quad & \sum_{i=1}^{N} f_i({\x}_i) + I_{\mathcal{Z}}({\z}),\nonumber \\
\text{s.t.}\quad &{\A}_i{\x}_i + {\z}_i  = \frac{\mathbf{c}}{N}, \quad i = 1,\ldots,N,
\end{align}
where $\z = ({\z}_1^{\top}, \ldots, {\z}_N^{\top})^{\top}$ is partitioned conformably according to ${\x}$, and $I_{\mathcal{Z}}({\z})$ is the indicator function of the convex set $\mathcal{Z}$, i.e.\ $I_{\mathcal{Z}}({\z}) = 0$ for ${\z} \in \mathcal{Z} = \{{\z}| \sum_{i=1}^{N}{\z}_i = 0\}$ and $I_{\mathcal{Z}}({\z}) = \infty$ otherwise. The augmented Lagrangian function is
\begin{align}
\mathcal{L}_{\rho} &= \sum_{i=1}^{N} f_i({\x}_i) + I_{\mathcal{Z}}({\z}) - \sum_{i=1}^{N}{\blambda}_{i}^{\top}({\A}_i{\x}_i + {\z}_i  - \frac{\mathbf{c}}{N}) \nonumber \\
& + \frac{\rho}{2}\sum_{i=1}^{N} \Vert {\A}_i{\x}_i + {\z}_i  - \frac{\mathbf{c}}{N} \Vert_2^2,
\end{align}
where we have two groups of variables, $\{{\x}_1, \ldots, {\x}_N\}$ and $\{{\z}_1, \ldots, {\z}_N\}$. Hence, we can apply the two-block ADMM to update these two groups of variables iteratively, i.e,  we can first update group $\{{\x}_i\}$ and then update group $\{{\z}_i\}$. In each group, ${\x}_i$ and ${\z}_i$ can be updated concurrently in parallel at each iteration. In particular, the update rules for ${\x}_i$ and ${\z}_i$ are
\begin{equation}
\left\{\begin{array}{ll}
{\x}_i^{k+1} = \argmin_{{\x}_i}\mathcal{L}_{\rho} ({\x}_i,{\z}_i^{k},{\blambda}_{i}^{k}),&\\
{\z}_i^{k+1} = \argmin_{{\z}_i}\mathcal{L}_{\rho} ({\x}_1^{k+1},{\z}_i,{\blambda}_{i}^{k}),& \forall i = 1,\ldots,N,\\
{\blambda}_{i}^{k+1} = {\blambda}_{i}^{k} - {\rho}({\A}_i{\x}_i + {\z}_i  - \frac{\mathbf{c}}{N}).&
\end{array} \right.
\end{equation}
The variable splitting ADMM is illustrated in Algorithm \ref{alg:A4}. The relationship between this splitting scheme and the Jacobian splitting scheme has been outlined in the following work \cite{BMX13}.  Algorithm \ref{alg:A4} enjoys the convergence rates of the 2-block ADMM. However, the number of variables and constraints will increase substantially when $N$ is large, which will impact the efficiency and incur significant burden for the computation.

\begin{algorithm}[t]
\caption{Variable Splitting Multi-block ADMM}\label{alg:A4}
\begin{algorithmic}
\STATE Initialize: ${\x}^{0}$, ${\z}^{0}$, ${\blambda}^{0}$, $\rho>0$;

\FOR{$k=0,1,\ldots$}

\FOR{$i=1,\ldots,N$}

\STATE \COMMENT{${\x}_i$, ${\z}_i$ and ${\blambda}_{i}$ are updated \textbf{concurrently}.}

\STATE ${\x}_i^{k+1} = \argmin_{{\x}_i}\mathcal{L}_{\rho} ({\x}_i,{\z}_i^{k},{\blambda}_{i}^{k})$;

\STATE ${\z}_i^{k+1} = \argmin_{{\z}_i}\mathcal{L}_{\rho} ({\x}_1^{k+1},{\z}_i,{\blambda}_{i}^{k})$;

\STATE ${\blambda}_{i}^{k+1} = {\blambda}_{i}^{k} - {\rho}({\A}_i{\x}_i + {\z}_i  - \frac{\mathbf{c}}{N})$;

\ENDFOR

\ENDFOR

\end{algorithmic}
\end{algorithm}

\subsection{ADMM with Gaussian Back Substitution}
\label{sec:GSADMM}
Many efforts have been made to improve the convergence of the Guass-Seidel type multi-block ADMM \cite{BXM12,MTXMSZ14}. In this part, we describe the ADMM with Gaussian back substitution \cite{BXM12}, which asserts that if a new iterate is generated by correcting the output of Algorithm \ref{alg:A2} with a Gaussian back substitution procedure, then the sequence of iterates converges to a solution of problem (\ref{eqn:intro}). We first define vector $\mathbf{v} =({\x}_2^{\top}, \ldots, {\x}_N^{\top}, {\blambda}^{\top})^{\top}$, vector $\tilde{\mathbf{v}} =(\tilde{\x}_2^{\top}, \ldots, \tilde{\x}_N^{\top}, \tilde{\blambda}^{\top})^{\top}$, matrix $\mathbf{H} = \text{diag}(\rho {\A}_2^{\top}{\A}_2,\ldots,\rho {\A}_N^{\top}{\A}_N,\frac{1}{\rho}\mathbf{I}_{m})$ and $\mathbf{M}$ as
\begin{equation}
\mathbf{M} =
\left( \begin{array}{ccccc}
\rho {\A}_2^{\top}{\A}_2 & 0 & \ldots & \ldots & 0 \\
\rho {\A}_3^{\top}{\A}_2 & \rho {\A}_3^{\top}{\A}_3 & \ddots & &\vdots \\
\vdots & \vdots & \ddots & \ddots & \vdots\\
\rho {\A}_N^{\top}{\A}_2 & \rho {\A}_N^{\top}{\A}_3 & \ldots & \rho {\A}_N^{\top}{\A}_N & 0\\
0&0&\ldots&0&\frac{1}{\rho}\mathbf{I}_{m}
\end{array} \right).
\end{equation}

Each iteration of the ADMM with Gaussian back substitution consists of two procedures: a prediction procedure and a correction procedure. The $\tilde{\mathbf{v}}$ is generated by Algorithm \ref{alg:A2}. In particular, $\tilde{\x}_i$ is updated sequentially as
\begin{equation}
\tilde{\x}_i^{k} = \argmin_{\tilde{\x}_i}\mathcal{L}_{\rho}(\{\tilde{\x}_j^{k}\}_{j<i},{\x}_i,\{{\x}_j^{k}\}_{j>i},{\blambda}^k),
\end{equation}
where the prediction procedure is performed in a forward manner, i.e. from the first to the last block and to the Lagrangian multiplier. Note that the newly generated $\tilde{\x}_i$ are used in the update of the next block in accordance with the Gauss-Seidel update fashion. After the update of the Lagrangian multiplier, the correction procedure is performed update $\mathbf{v}$ as
\begin{equation}
\mathbf{H}^{-1}\mathbf{M}^{\top}(\mathbf{v}^{k+1}-\mathbf{v}^{k}) = \alpha(\tilde{\mathbf{v}}^{k}-\mathbf{v}^{k}),
\end{equation}
where $\mathbf{H}^{-1}\mathbf{M}^{\top}$ is an upper-triangular block matrix according to the definition of $\mathbf{H}$ and $\mathbf{M}$. This implies that the update of correction procedure is in a backward fashion, i.e., first update the Lagrangian multiplier, and then update ${\x}_i$ from the last block to the first block sequentially. Note that an additional assumption that ${\A}_i^{\top}{\A}_i(i = 1,2, \ldots, N)$ are nonsingular are made here. ${\x}_1$ serves as an intermediate variable and is unchanged during the correction procedure. The algorithm is illustrated in Algorithm \ref{alg:A5}.

The global convergence of the ADMM with Gaussian back substitution is proved in \cite{BXM12}, and the convergence rate and iteration complexity are addressed in \cite{BMX12}.

\begin{algorithm}[t]
\caption{ADMM with Gaussian Back Substitution}\label{alg:A5}
\begin{algorithmic}
\STATE Initialize: ${\x}^{0}$, $\tilde{\x}^{0}$, ${\blambda}^{0}$, $\tilde{\blambda}^{0}$, $\rho>0$, $\alpha \in (0,1)$;

\FOR{$k=0,1,\ldots$}

\FOR{$i=1,\ldots,N$}

\STATE \COMMENT{${\x}_i$ is updated \textbf{sequentially}.}

\STATE $\tilde{\x}_i^{k} = \argmin_{\tilde{\x}_i}\mathcal{L}_{\rho}(\{\tilde{\x}_j^{k}\}_{j<i},{\x}_i,\{{\x}_j^{k}\}_{j>i},{\blambda}^k)$;

\ENDFOR

\STATE $\tilde{\blambda}^{k+1} = {\blambda}^{k} - {\rho}(\sum_{i=1}^{N} {\A}_i\tilde{\x}_i^{k+1} - \mathbf{c})$;

\STATE

\STATE \COMMENT{Gaussian back substitution correction step}

\STATE $\mathbf{H}^{-1}\mathbf{M}^{\top}(\mathbf{v}^{k+1}-\mathbf{v}^{k}) = \alpha(\tilde{\mathbf{v}}^{k}-\mathbf{v}^{k})$;

\STATE ${\x}_1^{k+1} = \tilde{\x}_1^{k}$;
\ENDFOR

\end{algorithmic}
\end{algorithm}

\subsection{Proximal Jacobian ADMM}
\label{sec:PJADMM}
The other type of modification on the ADMM for the multi-block setting is based on the Jacobian iteration scheme \cite{BLX13,WMZW14,HAZ14,BMX13}. Since the Guass-Seidel update is performed sequentially and is not amenable for parallelization,  Jacobian type iteration is preferred for distributed and parallel optimization.  In this subsection, we describe the proximal Jacobian ADMM \cite{WMZW14}, in which a proximal term \cite{NS13} is added in the update compare with that of Algorithm \ref{alg:A3} to improve convergence. In particular, the update of ${\x}_i$ is
\begin{equation}
{\x}_i^{k+1} \! = \! \argmin_{{\x}_i}\mathcal{L}_{\rho}({\x}_i,\{{\x}_j^{k}\}_{j \neq i},{\blambda}^k) \! + \! \frac{1}{2}\Vert {\x}_i \! - \! {\x}_i^{k}\Vert_{\mathbf{P}_i}^{2},
\end{equation}
where $\Vert {\x}_i \Vert_{\mathbf{P}_i}^{2} = {\x}_i^{\top}\mathbf{P}_i{\x}_i$ for some symmetric and positive semi-definite matrix $\mathbf{P}_i \succeq 0$. The involvement of the proximal term can make the subproblem of ${\x}_i$ strictly or strongly convex, and thus make the problem more stable. Moreover, multiple choice of $\mathbf{P}_i$ can make the subproblems easier to solve. The update of the Lagrangian multiplier is
\begin{equation}
{\blambda}^{k+1} = {\blambda}^{k} - {\gamma}{\rho}(\sum_{i=1}^{N} {\A}_i{\x}_i^{k+1} - \mathbf{c}),
\end{equation}
where $\gamma >0$ is the damping parameter and the algorithm is illustrate in Algorithm \ref{alg:A6}.

The global convergence of the proximal Jacobian ADMM which is proved in \cite{WMZW14}. Moreover, it enjoys a convergence rate of $o(1/k)$ under conditions on $\mathbf{P}_i$ and $\gamma$.

\begin{algorithm}[t]
\caption{Proximal Jacobian ADMM}\label{alg:A6}
\begin{algorithmic}
\STATE Initialize: ${\x}^{0}$, ${\blambda}^{0}$, $\rho>0$, $\gamma >0$;

\FOR{$k=0,1,\ldots$}

\FOR{$i=1,\ldots,N$}

\STATE \COMMENT{${\x}_i$ is updated \textbf{concurrently}.}

\STATE ${\x}_i^{k+1} \! = \! \argmin_{{\x}_i}\mathcal{L}_{\rho}({\x}_i,\{{\x}_j^{k}\}_{j \neq i},{\blambda}^k) \! + \! \frac{1}{2}\Vert {\x}_i \! - \! {\x}_i^{k}\Vert_{\mathbf{P}_i}^{2}$;

\ENDFOR

\STATE ${\blambda}^{k+1} = {\blambda}^{k} - {\gamma}{\rho}(\sum_{i=1}^{N} {\A}_i{\x}_i^{k+1} - \mathbf{c})$;

\ENDFOR
\end{algorithmic}
\end{algorithm}

\subsection{Implementations}
The recent development in high performance computing (HPC) and cloud computing paradigm provides a flexible and efficient solution for deploying the large-scale optimization algorithms. In this part, we describe possible implementation approaches of those distributed and parallel algorithms on current mainstream large scale computing facilities.

One possible implementation utilizes available computing incentive techniques and tools like MPI, OpenMP, and OpenCL. The MPI is a language-independent protocol used for inter-process communications on distributed memory computing platform, and is widely used for high-performance parallel computing today. The (multi-block) ADMM using MPI has been implemented in \cite{BPCPE10} and \cite{ZMW13}. Besides, the OpenMP, which is a shared memory multiprocessing parallel computing paradigm, and the OpenCL, which is a heterogenous distributed-shared memory parallel computing paradigm that incorporate CPUs and GPUs, also promise to implement distributed and parallel optimization algorithms on HPC. It is expected that supercomputers will reach one exaflops ($10^{18}$ FLOPS) and even zettaflops ($10^{21}$ FLOPS) in the near feature, which will largely enhance the computing capacity and significantly expedite program execution.

Another possible approach exploits the ease-of-use cloud computing engine like Hadoop MapReduce and Apache Spark. The amount of cloud infrastructure available for Hadoop MapReduce makes it convenient to use for large problems, though it is awkward to express ADMM in MapReduce since it is not designed for iterative tasks. Apache Spark's in-memory computing feature enables it to run iterative optimizations much faster than Hadoop, and is now prevalent for large-scale machine learning and optimization task on clusters\cite{MMMSI10}. This implementation approach is much simpler than previous computing incentive techniques and tools and promise to implementation of the large-scale distributed and parallel computation algorithms based on ADMM. The advances in the cloud/cluster computing engine provides a simple method to implement the large-scale data processing, and recently Google, Baidu and Alibaba are also developing and deploying massive cluster computing engines to perform the large-scale distributed and parallel computation.

\begin{figure}[t]
    \centering
    \includegraphics[width=0.48\textwidth]{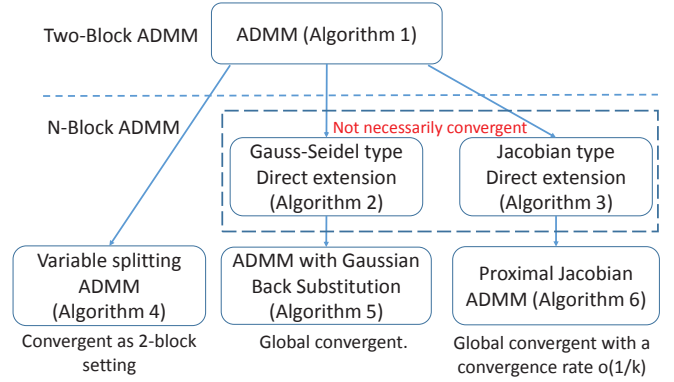}
    \caption{An illustration of the relationships among Algorithms.}
    \label{fig:alg}
\end{figure}

Now we have finished the review of distributed and parallel optimization methods based on ADMM, and we summarize the relationships among all Algorithms in Figure \ref{fig:alg}.

\section{Communication Network Applications}
\label{sec:application}
In this section, we review several applications of distributed and parallel optimization in communication networks. In particular, we describe the security constrained optimal power flow problem \cite{LAWZ13,SMERS14} in smart grids and the mobile data offloading in SDN \cite{LXMGZ15} based on ADMM.

\subsection{Security Constrained Optimal Power Flow}
\label{sec:SCOPF}
In this subsection, we consider the distributed and parallel approach for security constrained optimal power flow problem (SCOPF) \cite{LAWZ13, SMERS14}. The SCOPF is an extension of the conventional optimal power flow (OPF) problem, whose objective is to determine a generation schedule that minimizes the system operating cost while satisfying the system operation constraints such as hourly load demand, fuel limitations, environmental constraints and network security requirements.

\begin{figure}[t]
    \centering
    \includegraphics[width=0.48\textwidth]{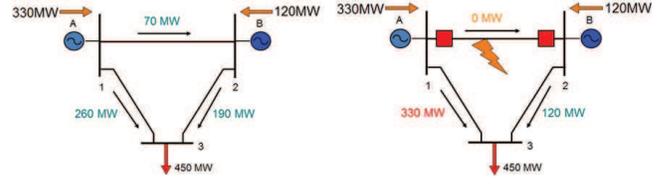}
    \caption{Example for the security constrained optimal power flow problem.}
    \label{fig:Contigency_example}
\end{figure}

An illustrative example of SCOPF is shown in Figure \ref{fig:Contigency_example}. There are $3$ buses with limit $300$MW, $2$ generators ($330$MW and $120$MW) and a load $450$MW in the system. In the left figure, it is an example for traditional optimal power flow problem (OPF) without considering the security constraint. If the line between buses $1$ and $2$ breaks, the line between buses $1$ and $3$ cannot afford $330$MW ($ \le 300$MW), and consequently it breaks. Then generator $B$ cannot afford the load, and as a result line between buses $2$ and $3$ breaks. We can see from this example why we need to consider the security constraint so as to avoid large area blackouts.

In \cite{LAWZ13}, the general form of SCOPF can be formulated as follows
\begin{align}
\label{chapter3:eqn:SCOPF}\min_{\mathbf{x}^0,\ldots,\mathbf{x}^C;\mathbf{u}^0,\ldots,\mathbf{u}^C} & \quad f^0(\mathbf{x}^0,\mathbf{u}^0) \\
\text{s.t.} & \label{chapter3:eqn:pre1} \quad \mathbf{g}^0(\mathbf{x}^0,\mathbf{u}^0) = 0,\\
                  & \label{chapter3:eqn:pre2} \quad \mathbf{h}^0(\mathbf{x}^0,\mathbf{u}^0) \le 0,\\
                  & \label{chapter3:eqn:pos1} \quad \mathbf{g}^c(\mathbf{x}^c,\mathbf{u}^c) = 0,\\
                  & \label{chapter3:eqn:pos2} \quad \mathbf{h}^c(\mathbf{x}^c,\mathbf{u}^c) \le 0, \text{ and}\\
                  & \label{chapter3:eqn:pos3} \quad \vert \mathbf{u}^0 - \mathbf{u}^c \vert \le \mathbf{\Delta}_c, \quad  c = 1,\ldots,C,
\end{align}
where $f^0$ is the objective function, which (\ref{chapter3:eqn:SCOPF}) aims to maximize the total social welfare or equivalently minimize offer-based energy and production cost, $\mathbf{x}^c$ is the vector of state variables, which includes voltage magnitudes and angles at all buses, and $\mathbf{u}^c$ is the vector of control variables, which can be generator real powers or terminal voltages. The superscript $c = 0$ corresponds to the pre-contingency configuration, and $c = 1,\ldots,C$ correspond to different post-contingency configurations. In addition, $\mathbf{\Delta}_c$ is the maximum allowed adjustment between the normal and contingency states for contingency $c$.

In the conventional SCOPF problem, the equality constraints \ref{chapter3:eqn:pos1} on $\mathbf{g}^c, c = 0,\ldots,C$, represent the system nodal power flow balance over the entire grid, and the inequality constraints \ref{chapter3:eqn:pos2} on $\mathbf{h}^c, c = 0,\ldots,C$, represent the physical limits on the equipment, such as the operational limits on the branch currents and bounds on the generator power outputs. Constraints (\ref{chapter3:eqn:pre1})-(\ref{chapter3:eqn:pre2}) capture the economic dispatch and  enforce the feasibility of the pre-contingency state. Constraints (\ref{chapter3:eqn:pos1})-(\ref{chapter3:eqn:pos2}) incorporate the security-constrained dispatch and enforce the feasibility of the post-contingency state. Constraint (\ref{chapter3:eqn:pos3}) introduces the security-constrained dispatch with rescheduling, which couples control variables of pre-contingency and post-contingency states and prevents unrealistic post-contingency corrective actions. Note that there are some variations on the objective function and constraints of the SCOPF problem, and we focus on the above conventional formulation in this subsection.

Following the standard approach to formulating the SCOPF problem, the objective here is to minimize the cost of generation while safeguarding the power system sustainability. For the sake of simplicity and computational tractability, constraints (\ref{chapter3:eqn:pre1})-(\ref{chapter3:eqn:pos2}) are modeled with the linear DC load flow, and we assume that the list of contingencies is given. Thus, assuming a DC power network modeling and neglecting all shunt elements, the standard SCOPF problem can be simplified to the following optimization problem
\begin{align}
\label{chapter3:eqn:DCSCOPF}\min_{\boldsymbol \theta^0,\ldots,\boldsymbol \theta^C;\mathbf{P}^{g,0},\ldots,\mathbf{P}^{g,C}} & \quad \sum_{i \in \mathcal{G}}f_i^g(\mathbf{P}_i^{g,0}) \\
\text{s.t.} & \label{chapter3:eqn:DCpre1} \quad \mathbf{B}_{bus}^0 \boldsymbol \theta^0 + \mathbf{P}^{d,0} - \mathbf{A}^{g,0}\mathbf{P}^{g,0}= 0,\\
                  & \label{chapter3:eqn:DCpos1} \quad \mathbf{B}_{bus}^c \boldsymbol \theta^c + \mathbf{P}^{d,c} - \mathbf{A}^{g,c}\mathbf{P}^{g,c}= 0,\\
                  & \label{chapter3:eqn:DCpre2} \quad \vert \mathbf{B}_{f}^0 \boldsymbol \theta^0 \vert - \mathbf{F}_{max} \le 0,\\
                  & \label{chapter3:eqn:DCpos2} \quad \vert \mathbf{B}_{f}^c \boldsymbol \theta^c\vert - \mathbf{F}_{max} \le 0,\\
                  & \label{chapter3:eqn:DCpre3} \quad \underline{\mathbf{P}^{g,0}} \le \mathbf{P}^{g,0} \le \overline{\mathbf{P}^{g,0}},\\
                  & \label{chapter3:eqn:DCpos3} \quad \underline{\mathbf{P}^{g,c}} \le \mathbf{P}^{g,c} \le \overline{\mathbf{P}^{g,c}},\\
                  & \label{chapter3:eqn:DCpos4} \quad \vert \mathbf{P}^{g,0} - \mathbf{P}^{g,c} \vert \le \mathbf{\Delta}_c, \text{ and}\\
                  & \quad i \in \mathcal{G}, \quad c = 1,\ldots,C,
\end{align}
where the notation is given in Table \ref{chapter3:table:Notation}.

The solution to (\ref{chapter3:eqn:DCSCOPF}) ensures economical dispatch while guaranteing power system security, by taking into account a set of postulated contingencies. The major challenge of SCOPF is the problem size, especially for large systems with numerous contingency cases to be considered. Directly solving the SCOPF problem by simultaneously imposing all post-contingency constraints will result in prohibitive memory requirements and a substantial CPU burden. The proposed distributed optimization method is based on the ADMM. However, the optimization problem (\ref{chapter3:eqn:DCSCOPF}) cannot be readily solved using ADMM, since the constraint (\ref{chapter3:eqn:DCpos4}) couples the pre-contingency and post-contingency variables, and the inequalities make the problem even more complicated.  \begin{table}
\centering
\caption{Notation definitions.}\label{chapter3:table:Notation}
\begin{tabular}{|l|p{5.5cm}|}
\hline
$\mathcal{G}$ & Set of generators\\
$\mathcal{N}$ & Set of buses\\
$\mathcal{B}$ & Set of branches\\
$\boldsymbol \theta^c \in \mathbb{R}^{\vert \mathcal{N} \vert}$ & Vector of voltage angles\\
$\mathbf{P}^{g,c} \in \mathbb{R}^{\vert \mathcal{G} \vert}$ & Vector of real power flows\\
$f_i^g$ & Generation cost function\\
$\mathbf{P}_i^{g,0}$ & Displaceable real power of each individual generation unit $i$ for the pre-contingency configuration\\
$\mathbf{B}_{bus}^c \in \mathbb{R}^{\vert \mathcal{N} \vert \times \vert \mathcal{N} \vert}$ & Power network system admittance matrix\\
$\mathbf{B}_{f}^c \in \mathbb{R}^{\vert \mathcal{B} \vert \times \vert \mathcal{N} \vert}$ & Branch admittance matrix\\
$\mathbf{P}^{d,c} \in \mathbb{R}^{\vert \mathcal{N} \vert}$ & Real power demand\\
$\mathbf{A}^{g,c} \in \mathbb{R}^{\vert \mathcal{N} \vert \times \vert \mathcal{G} \vert}$ & Sparse generator connection matrix, whose $(i,j)$-th element is 1 if generator $j$ is located at bus $i$ and 0 otherwise\\
$\mathbf{F}_{max}$ & Vector for the maximum power flow\\
$\overline{\mathbf{P}^{g,c}}$ & Upper bound on real power generation\\
$\underline{\mathbf{P}^{g,c}}$ & Lower bound on real power generation\\
$\mathbf{\Delta}_c$ & Pre-defined maximum allowed variation of power outputs\\
\hline
\end{tabular}
\end{table}To address these challenges, the optimization problem (\ref{chapter3:eqn:DCSCOPF}) can then be reformulated by introducing a slack variable $\mathbf{p}^c \in \mathbb{R}^{\vert \mathcal{G} \vert}$
\begin{align}
\text{minimize}   & \quad (\ref{chapter3:eqn:DCSCOPF})\\
\text{subject to} & \quad \text{Constraints (\ref{chapter3:eqn:DCpre1})-(\ref{chapter3:eqn:DCpos3}}),\\
                  & \label{chapter3:eqn:coupling}\quad \mathbf{P}^{g,0} - \mathbf{P}^{g,c} + \mathbf{p}^{c} = \mathbf{\Delta}_c, \text{ and}\\
                  & \label{chapter3:eqn:slack} \quad 0 \le \mathbf{p}^{c} \le 2\mathbf{\Delta}_c,\quad c = 1,\ldots,C.
\end{align}

The above optimization problem can be solved distributively using ADMM. The scaled augmented Lagrangian can be calculated as
\begin{align}
& \mathcal{L}_{\rho} (\{\mathbf{P}^{g,c}\}_{c=1}^{C};\{\mathbf{p}^c\}_{c=1}^{C}; \{\boldsymbol \mu^c\}_{c=1}^{C}) = \nonumber \\
& \sum_{i \in \mathcal{G}}f_i^g(\mathbf{P}_i^{g,0})\! + \! \sum_{c = 1}^{C}\frac{\rho^c}{2}\Vert \mathbf{P}^{g,0} \! - \! \mathbf{P}^{g,c} \!+ \! \mathbf{p}^{c} \!-\! \mathbf{\Delta}_c \!+\! \boldsymbol \mu^c \Vert_2^2.
\end{align}

\begin{figure*}[t]
  \centering
  \subfigure[Computing time for the IEEE 57 bus system]{
    \label{chapter3:fig:Case57RT} 
    \includegraphics[height = 1.8in, width=0.32\linewidth]{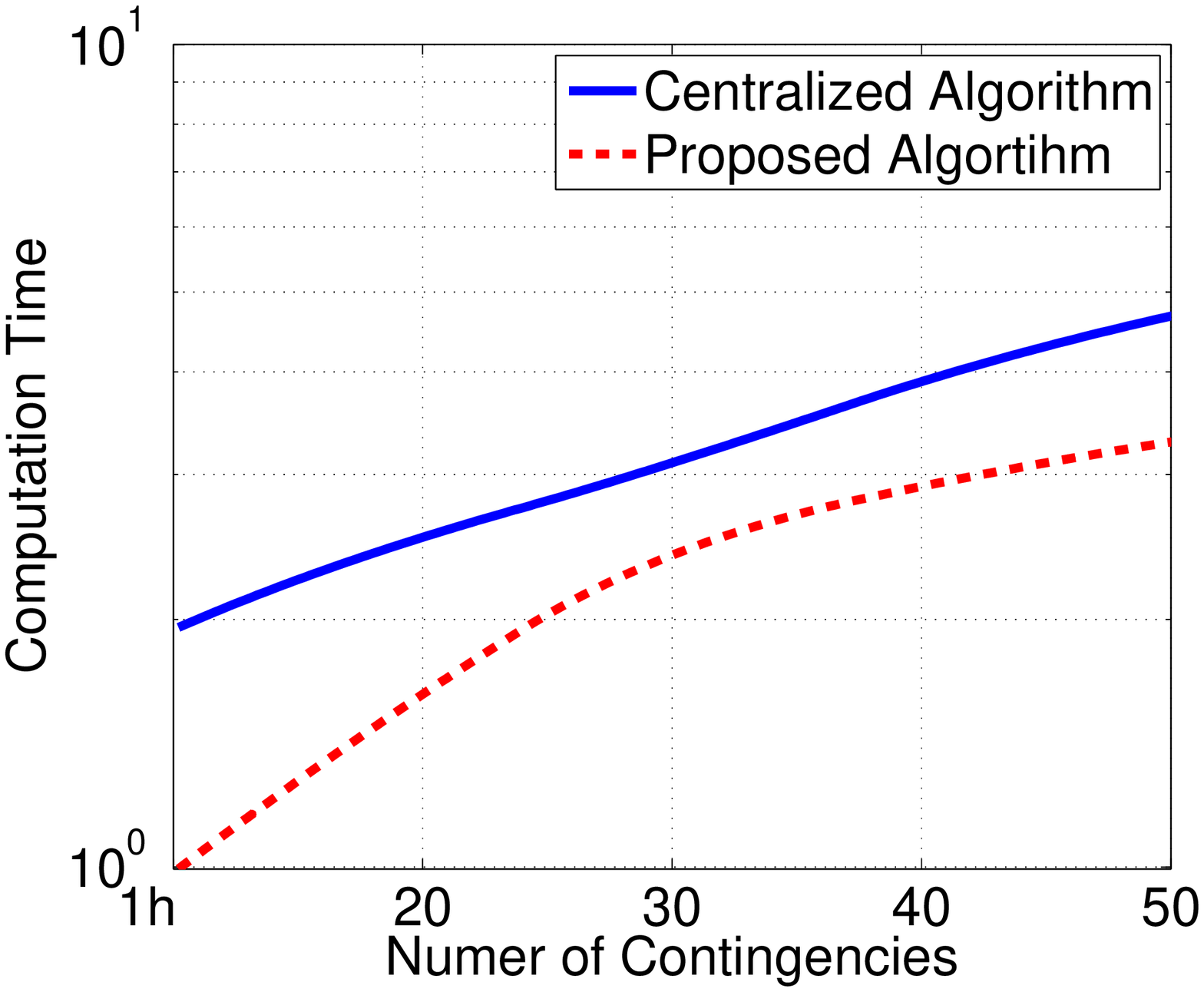}}
  \subfigure[Computing time for the IEEE 118 bus system]{
    \label{chapter3:fig:Case118RT} 
    \includegraphics[height = 1.8in, width= 0.32\linewidth]{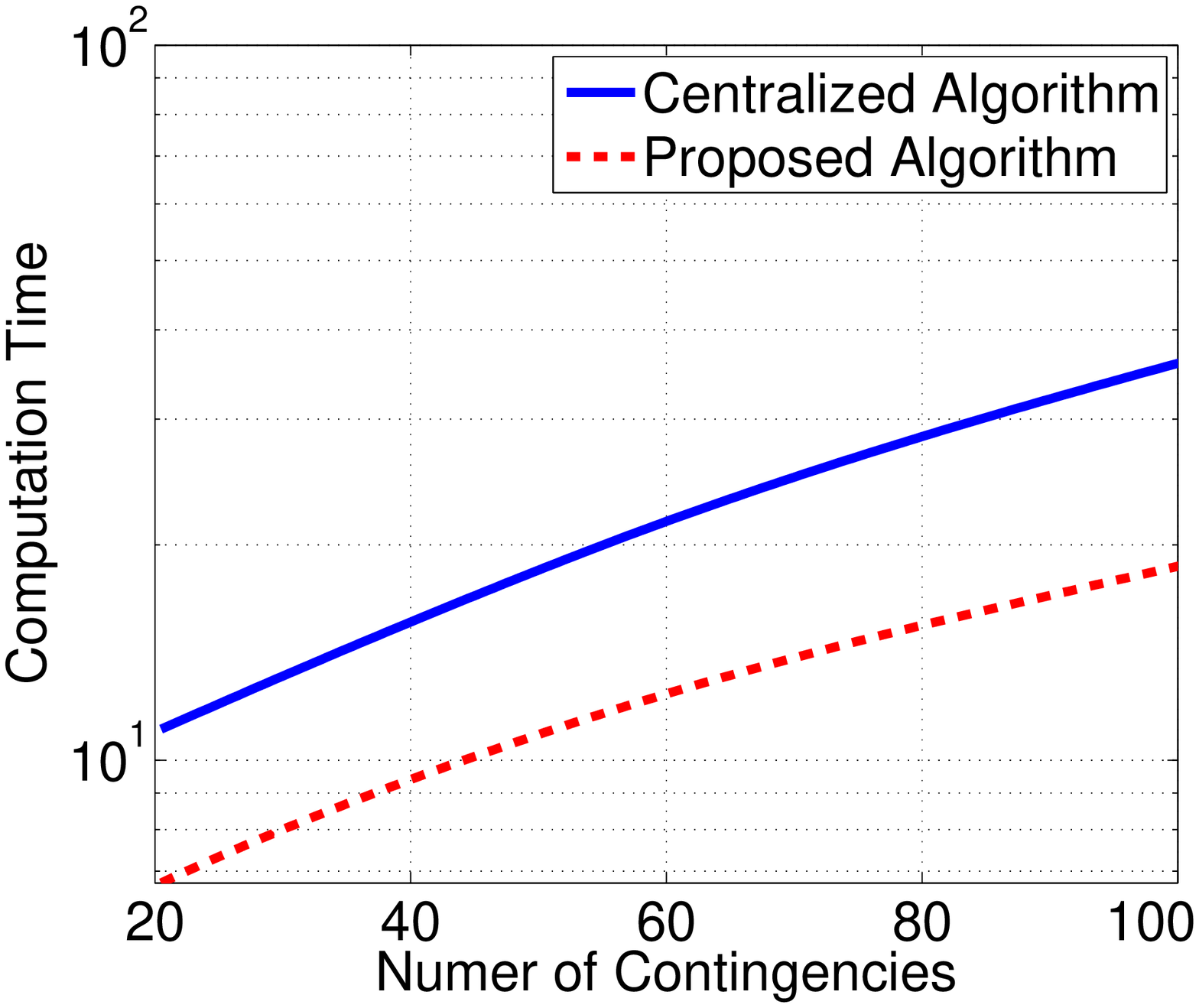}}
  \subfigure[Computing time for the IEEE 300 bus system]{
    \label{chapter3:fig:Case300T} 
    \includegraphics[height = 1.8in, width=0.32\linewidth]{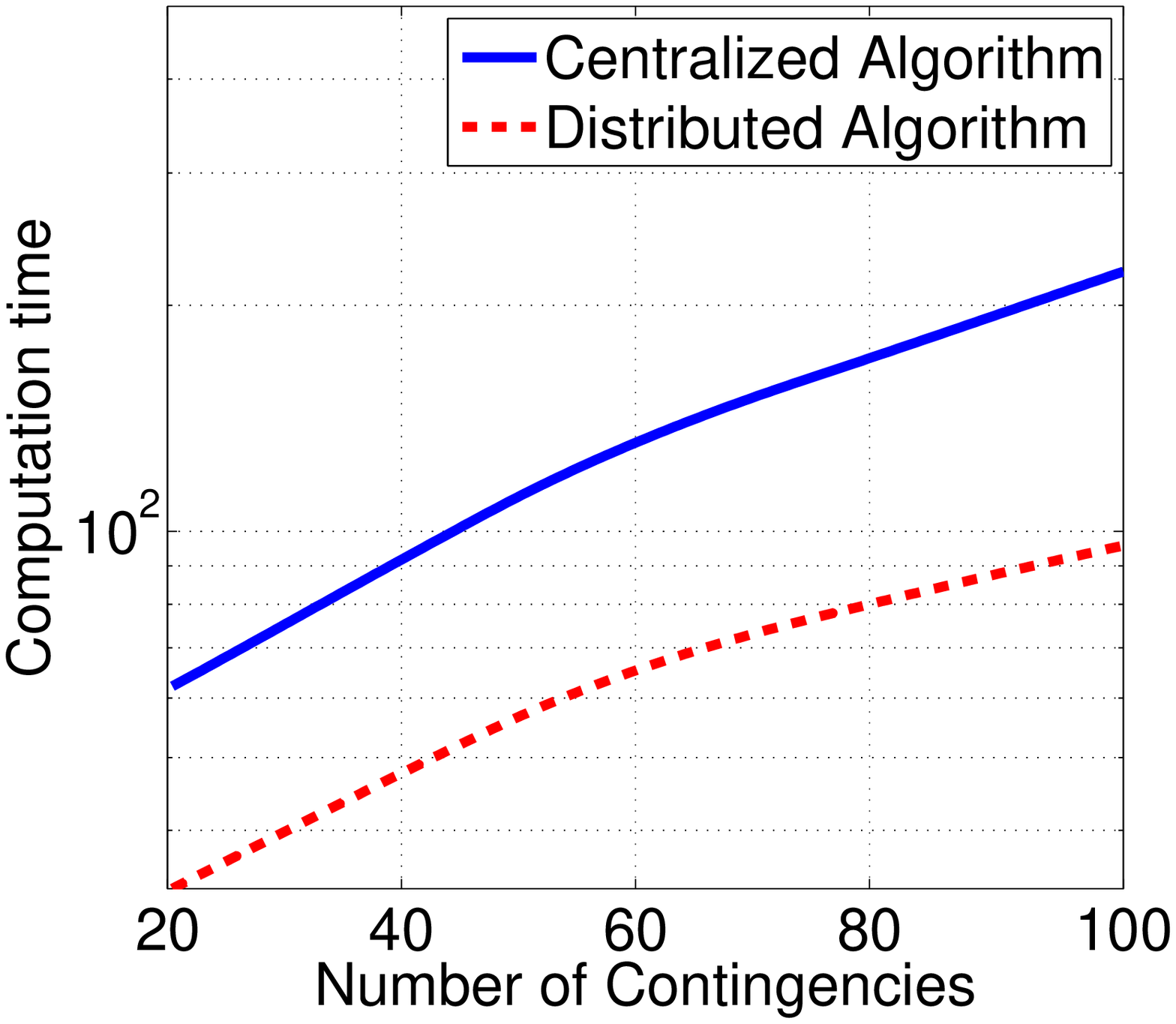}}
  \caption{Computing time for the IEEE 57 bus system, IEEE 118 bus system and IEEE 300 bus system with different numbers of contingency cases.}
  \label{chapter3:fig:CompT} 
\end{figure*}

The optimization variables $\mathbf{P}^{g,0}, \mathbf{P}^{g,c}$, and $\mathbf{p}^{c}$ are arranged into two groups, $\{ \mathbf{P}^{g,0} \}$ and $\{ \mathbf{P}^{g,c}, \mathbf{p}^{c}\}$, and updated iteratively. The variables in each group are optimized in parallel on distributed computing nodes, and coordinated by the dual variable vector $\boldsymbol \mu^c$ during each iteration.

At the $k^{th}$ iteration, the $\mathbf{P}^{g,0}$-update solves the base scenario with squared regularization terms enforced by the coupling constraints and expressed as
\begin{align}
&\mathbf{P}^{g,0}[k+1] = \argmin_{\mathbf{P}^{g,0}} \sum_{i \in \mathcal{G}}f_i^g(\mathbf{P}_i^{g,0}) \nonumber \\
& + \sum_{c = 1}^{C}\frac{\rho^c}{2}\Vert \mathbf{P}^{g,0} - \mathbf{P}^{g,c}[k] + \mathbf{p}^{c}[k] - \mathbf{\Delta}_c + \boldsymbol \mu^c[k] \Vert_2^2, \nonumber \\
&\text{subject to} \quad \text{Constraints} (\ref{chapter3:eqn:DCpre1}),(\ref{chapter3:eqn:DCpre2}), \text{and } (\ref{chapter3:eqn:DCpre3}).
\end{align}
The $\mathbf{P}^{g,c}$-updating solves a number of independent optimization subproblems correspond to post-contingency scenarios and can be calculated distributively at the $c^{th}$ computing nodes via
\begin{align}
& \label{chapter3:eqn:sub} \mathbf{P}^{g,c}[k+1] = \nonumber \\
& \argmin_{\mathbf{P}^{g,c},\mathbf{p}^{c}} \frac{\rho^c}{2}\Vert \mathbf{P}^{g,0}[k+1] - \mathbf{P}^{g,c} + \mathbf{p}^{c} - \mathbf{\Delta}_c + \boldsymbol \mu^c[k] \Vert_2^2, \nonumber \\
&\text{subject to} \quad \text{Constraints} (\ref{chapter3:eqn:DCpos1}),(\ref{chapter3:eqn:DCpos2}),(\ref{chapter3:eqn:DCpos3}), \text{and } (\ref{chapter3:eqn:slack}),
\end{align}
where the scaled dual variable vector is also updated locally at the $c^{th}$ computing utility as
\begin{equation}
\boldsymbol \mu^c[k+1] = \boldsymbol \mu^c[k] + \mathbf{P}^{g,0}[k+1] - \mathbf{P}^{g,c}[k+1] + \mathbf{p}^{c}[k+1] - \mathbf{\Delta}_c.
\end{equation}

At the $k^{th}$ iteration, the original problem is divided into $C+1$ subproblems of approximately the same size. The computing node handling $\mathbf{P}^{g,0}$ needs to communicate with all computing nodes solving (\ref{chapter3:eqn:sub}) during the iterations. The results of the $\mathbf{P}^{g,0}$-update, $\{\mathbf{P}^{g,0}\}$, will be distributed among the computing nodes for the $\mathbf{P}^{g,c}$-update. After the $\mathbf{P}^{g,c}$-update, the computed $\{\mathbf{P}^{g,c},\mathbf{p}^c, \boldsymbol \mu^c\}$ will be collected to calculate the pre-contingency control variables. The subproblem data are iteratively updated such the block-coupling constraints (\ref{chapter3:eqn:coupling}) are satisfied at the end. Note that since each of the subproblems is a smaller-scale OPF problem, existing techniques for OPF can be applied with minor modifications.The proposed algorithm is illustrated in Algorithm \ref{chapter3:alg:A3}.

Numerical studies are examined to evaluate the performance of the proposed algorithm. Three classical test systems are used: the IEEE 57 bus, the IEEE 118 bus, and the IEEE 300 bus. The computing time for test systems with different numbers of contingency cases is investigated and results are given in Figure \ref{chapter3:fig:CompT}. The number of contingencies is increased by $20\%$ each time and the computing time is recorded. It can be seen from these figures that with an increase in the number of contingency cases for the SCOPF problem, the computing time of the centralized algorithm increases much faster than that of the proposed algorithm. Thus, the proposed distributed algorithm is more scalable and stable than the centralized approach.

\begin{algorithm}[t]
\caption{Distributed SCOPF.}\label{chapter3:alg:A3}
\begin{algorithmic}
\STATE Input: $\mathbf{B}_{bus}^c$, $\mathbf{B}_{f}^c$, $\mathbf{A}^{g,c}$, $\mathbf{P}^{d,c}$, $\overline{\mathbf{P}^{g,c}}$, $\underline{\mathbf{P}^{g,c}}$, $\mathbf{\Delta}_c$;

\STATE Initialize: $\boldsymbol \theta^c$, $\mathbf{P}^{g,c}$, $\mathbf{p}^{c}$, $\boldsymbol \mu^c$, $\rho^c$, $k = 0$;

\WHILE {not converge}

\STATE {$\mathbf{P}^{g,0}$-update}:

\STATE $\mathbf{P}^{g,0}[k+1] = \argmin_{\mathbf{P}^{g,0}} \sum_{i \in \mathcal{G}}f_i^g(\mathbf{P}_i^{g,0})$

\STATE $+ \sum_{c = 1}^{C}\frac{\rho^c}{2}\Vert \mathbf{P}^{g,0} - \mathbf{P}^{g,c}[k] + \mathbf{p}^{c}[k] - \mathbf{\Delta}_c + \boldsymbol \mu^c[k] \Vert_2^2$

\STATE {subject to Constraints} (\ref{chapter3:eqn:DCpre1}),(\ref{chapter3:eqn:DCpre2}), and (\ref{chapter3:eqn:DCpre3}).

\STATE

\STATE {$\mathbf{P}^{g,c}$-update, distributively at each computing node}:

\STATE $\mathbf{P}^{g,c}[k+1] = \argmin_{\mathbf{P}^{g,c},\mathbf{p}^c} \frac{\rho^c}{2}\Vert \mathbf{P}^{g,0}[k+1] - \mathbf{P}^{g,c} + \mathbf{p}^{c} - \mathbf{\Delta}_c + \boldsymbol \mu^c[k] \Vert_2^2$

\STATE {subject to Constraints} (\ref{chapter3:eqn:DCpos1}),(\ref{chapter3:eqn:DCpos2}),(\ref{chapter3:eqn:DCpos3}), and (\ref{chapter3:eqn:slack}),

\STATE $\boldsymbol \mu^c[k+1] = \boldsymbol \mu^c[k] + \mathbf{P}^{g,0}[k+1] - \mathbf{P}^{g,c}[k+1] + \mathbf{p}^{c}[k+1] - \mathbf{\Delta}_c.$

\STATE {Adjust penalty parameter $\rho^c$ is necessary;}

\STATE $k = k + 1$;

\ENDWHILE

\RETURN $\boldsymbol \theta^c$, $\mathbf{P}^{g,c}$;

\STATE Output $\boldsymbol \theta^c$, $\mathbf{P}^{g,c}$;

\end{algorithmic}
\end{algorithm}

\subsection{Mobile Data Offloading in SDN}
We consider a mobile network which consists of $B$ cellular base stations (BSs) and $A$ access points (APs). A BS $b \in \{ 1, \ldots, B\}$ serves a group of mobile users and has the demand to offload its traffic to APs. An AP $a \in \{1, \ldots, A\}$ is a WiFi or femtocell AP which operates in a different frequency band and supply its bandwidth for data offloading. The maximum available capacity for data offloading of each AP $a$ is denoted by $C_a$. The SDN controller manages the BSs and APs through the access network discovery and selection function (ANDSF), and makes the mobile data offloading decisions according to various trigger criteria. Such criteria can be the number of mobile users per BS, available bandwidth/IP address of each BS, or aggregate number of flows on a specific port at a BS.

Let $\x_b = [x_{b1}, \ldots, x_{bA}]^{\top}$ represent the offloaded traffic of BS $b$, where $x_{ba}$ denotes the data of BS $b$ offloaded through AP $a$. Correspondingly, $\y_a = [y_{a1}, \ldots, y_{aB}]^{\top}$ represents the admitted traffic of AP $a$, where $y_{ab}$ represents the admitted data traffic from BS $b$. Generally, a feasible mobile data offloading decision exists when BSs and APs reach an agreement on the amount of offloading data, i.e., $x_{ba} = y_{ab}, \forall a$ and $\forall b$. We assume that the mobile data of BSs can be offloaded to all of the APs without loss of generality. Moreover, we assume that the time is slotted and during each slot duration the offloading demand from BSs is fixed. The SDN controller needs to find a feasible offloading schedule at the beginning of each time slot, while maximizing the utility of BSs at a reasonable cost of APs.

We denote BS $b$'s utility of offloading its traffic to APs by $U_b(\x_b)$, where $U_b(\cdot)$ is designed to be a non-decreasing, non-negative and concave function in $\x_b, \forall b$. For example, the function can be logarithmic, and the concavity is justified because of diminishing returns of the resources allocated to the offload data. Likewise, we use function $L_a(\y_a)$ to describe the AP $a$'s cost of helping BSs offload data, where $L_a(\cdot)$ is a non-decreasing, non-negative and convex function in $\y_a, \forall a$. The cost function can be a linear cost function, which means the total cost of APs will increase as the amount of admitted mobile data increases.

For the SDN controller, the total revenue for mobile data offloading is expressed as $\sum_{b=1}^{B}U_b({\x}_b) - \sum_{a=1}^{A}L_a({\y}_a)$. To maximize the total revenue, the equivalent minimization optimization problem can be formulated as,
\begin{align}
\label{eqn:NUM}
\min_{\{\x_1,\ldots,\x_B\},\{\y_1,\ldots,\y_A\}} & \quad \sum_{a=1}^{A}L_a({\y}_a) - \sum_{b=1}^{B}U_b({\x}_b), \\
\text{s.t} & \label{eqn:Capacity} \quad \sum_{b=1}^{B} y_{ab} \le C_a, \quad \forall a, \\
           & \label{eqn:Work} \quad x_{ba} = y_{ab}, \quad \forall a,b,
\end{align}
where (\ref{eqn:Capacity}) stands for the capacity constraint at each AP, and (\ref{eqn:Work}) represents the consensus of BSs and APs on the amount of mobile data.

\begin{figure}[t]
    \centering
    \includegraphics[height = 6cm,width=0.48\textwidth]{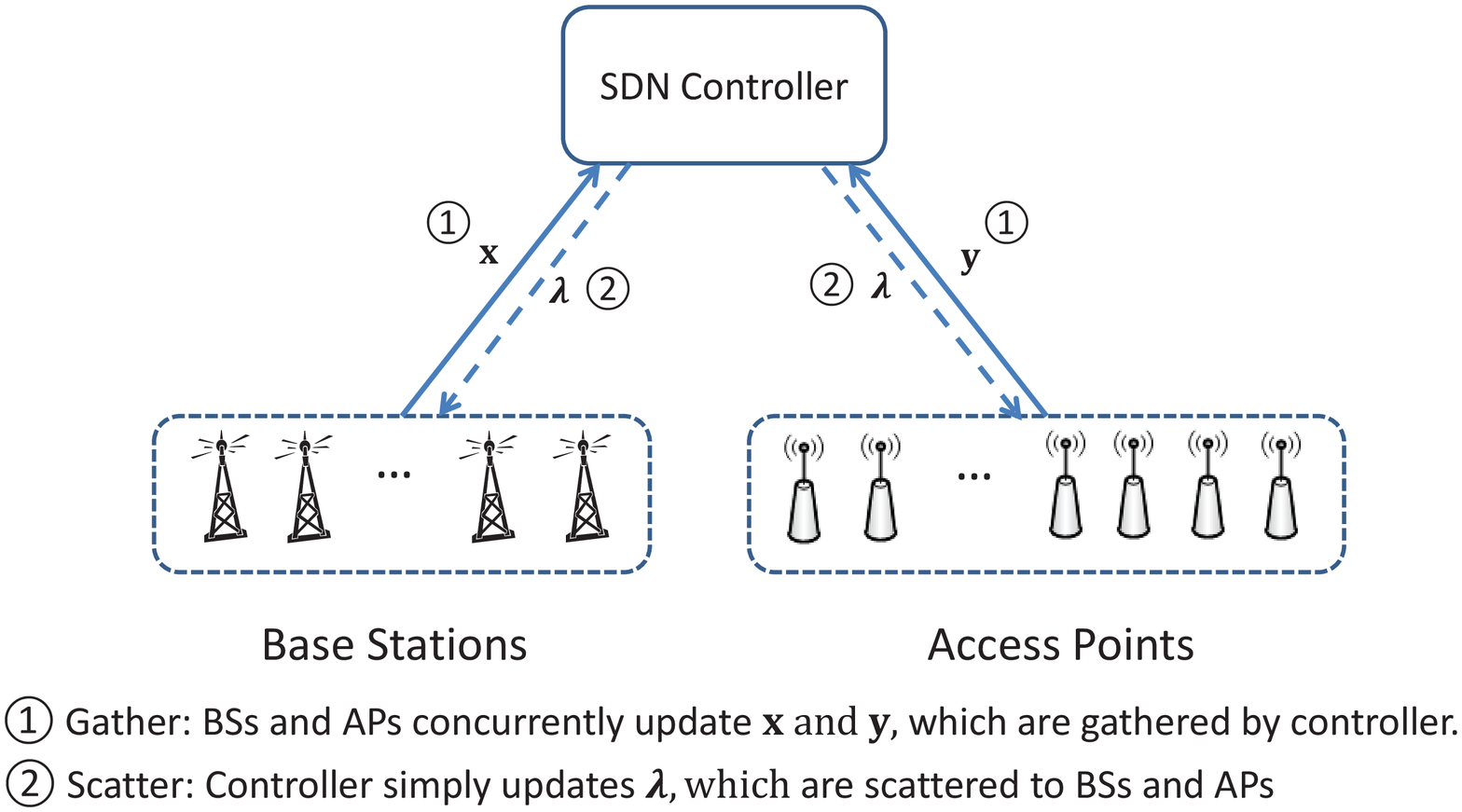}
    \caption{Distributed computing paradigm of proposed algorithm.}
    \label{fig:compu}
\end{figure}

We propose a fully distributed algorithm to solve the optimization problem (\ref{eqn:NUM}). The computing paradigm of the proposed algorithm is shown in Figure \ref{fig:compu} and can be described as follows. During each iteration, the BSs and APs update ${\x}$ and ${\y}$ concurrently. The updated ${\x}$ and ${\y}$ are gathered by the SDN controller, which performs a simple update on ${\blambda}$ and scatters the dual variables back to the BSs and APs. The iteration goes on until a consensus on the offloading demand and supply is reached. Specifically, we fist calculate the partial Lagrangian of (\ref{eqn:NUM}), which introduces the Lagrange multipliers only for constraint (\ref{eqn:Work}),
\begin{align}
&\mathcal{L}_{\rho} (\x,\y,{\blambda}) =  \sum_{a=1}^{A}L_a({\y}_a) - \sum_{b=1}^{B}U_b({\x}_b) \nonumber \\
&- \sum_{a = 1}^{A}\sum_{b=1}^{B} {\lambda}_{ab}(x_{ba}-y_{ab}) + \frac{\rho}{2}\sum_{a = 1}^{A}\sum_{b=1}^{B} \Vert x_{ba}-y_{ab}\Vert_2^{2},
\end{align}
where ${\blambda} \in \mathbb{R}^{AB}$ is the Lagrange multiplier and $\rho$ is the penalty parameter. The updates of BSs and APs can be performed concurrently according to the proximal Jacobian multi-block ADMM. We describe the update procedure of the BSs, APs and SDN controller as follows.

\textbf{Base Station Update}: At each BS $b$, the update rule can be expressed as,
\begin{align}
\label{eqn:BS}
{\x}_b^{k+1} \! = \! \argmin_{{\x}_b}( \! - \! U_b({\x}_b) \! +\! \frac{\rho}{2}\sum_{a=1}^{A}\Vert x_{ba}\!-\!p_{ab}^{k}\Vert_2^{2} \!+\!\frac{1}{2}\Vert{\x}_b \!-\! {\x}_b^{k}\Vert_{\mathbf{P}_i}^{2}),
\end{align}
where $\mathbf{P}_i = 0.1 \mathbf{I}$ and $\mathbf{I}$ is the identity matrix, and $p_{ab}^{k} = (y_{ab}^{k}+ \frac{{\lambda}_{ab}^{k}}{\rho}), \forall a$ is the `signal' sent from the SDN controller to BS $b$. The update (\ref{eqn:BS}) is a small scale unconstrained convex optimization problem. For each round of the update, it sends ${\x}_b$ of size $A$ to the SDN controller. Note that the update of each BS $b$ is performed independently and can be calculated locally. Once ${\x}_b$ is updated, it is sent to the SDN controller while the utility function $U_b(\cdot)$ is kept confidential.

\textbf{Access Point Update}: The update rule at each AP $a$ can be expressed as,
\begin{align}
\label{eqn:AP}
{\y}_a^{k+1}  &=  \argmin_{{\y}_b}(L_a({\y}_a)+\frac{\rho}{2}\sum_{b=1}^{B}\Vert y_{ab}- q_{ba}^{k}\Vert_2^{2} \nonumber \\
&+\frac{1}{2}\Vert{\y}_a - {\y}_a^{k}\Vert_{\mathbf{P}_i}^{2}), \quad \text{s.t.} \quad \sum_{b=1}^{B} y_{ab} \le C_a,
\end{align}
where $\mathbf{P}_i = 0.1 \mathbf{I}$ and $q_{ba}^{k} = (x_{ba}^{k} - \frac{{\lambda}_{ab}^{k}}{\rho}), \forall b$. $q_{ba}$ is the `signal' from the SDN controller to AP $a$. The update (\ref{eqn:AP}) is a small scale convex optimization problem with linear inequality constraints. For each round of the update, it sends ${\y}_a$ of size $B$ to the SDN controller. The update of ${\y}$ is also performed independently at each AP. During the update, the information of cost function $L_a(\cdot)$ is kept private. ${\y}_a$ is sent to the SDN controller once updated.

\textbf{SDN Controller Update}: At the SDN controller, the update rule can be expressed as,
\begin{equation}
\label{eqn:Con}
{\lambda}_{ab}^{k+1} = {\lambda}_{ab}^{k}-\gamma\rho\sum_{b=1}^{B}\sum_{a=1}^{A}(x_{ba}^{k+1}-y_{ab}^{k+1}).
\end{equation}
After gathering ${\x}$ and ${\y}$ from the BSs and APs, the SDN controller performs a simple update on the dual variable ${\blambda}$ by a simple algebra operation. After that, the `signal' variables $p_{ba}$ and $q_{ba}$ are scattered back to the corresponding BSs and APs, respectively. For each round of the update, it sends $p_{ba}, \forall a$ to each BS b, which is of size $A$, and sends $q_{ba}, \forall b$ to each AP a, which is of size $B$.

\begin{figure}[t]
    \centering
    \includegraphics[height = 6.4cm, width=0.48\textwidth]{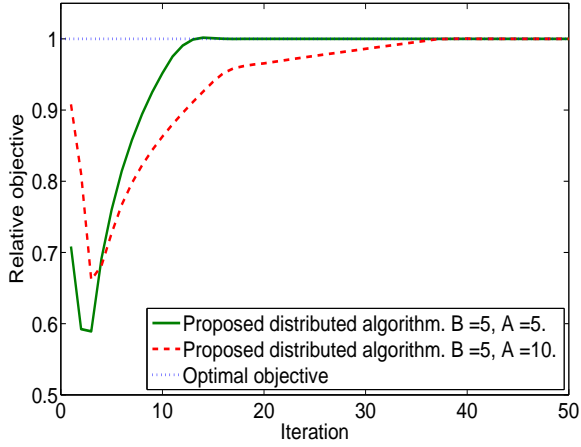}
    \caption{Convergence performance of the proposed algorithm by objective value when $(B = 5, A = 5)$ and $(B = 5, A = 10)$.}
    \label{fig:convergence}
\end{figure}

\textbf{Remark}  that in the Jacobian type update, the iterations of the BSs and APs are performed concurrently instead of consecutively in the Gauss-Seidel type update. There is no direct communication between the BSs and APs, which kept the intermediated update results of ${\x}$ and ${\y}$ confidential to each other. The updates at iteration $k+1$ only depends on its previous value at iteration $k$, which enables a fully distributed implementation.

At each iteration, the update operations at BSs and APs are quite simple. The update at each BS $b$ and AP $a$ are simple small scale convex optimization problems, which can be quickly solved by many off-the-shelf tools like CVX \cite{MS14}. As for the communication overhead, for each iteration the signaling between each BS and SDN controller is of the size $2A$ (size of ${\x}_b$ and $p_{ba}, \forall a$). Likewise, the signaling between each AP and SDN controller is of the size $2B$ (size of ${\y}_a$ and $q_{ba}, \forall b$). The sizes of those signaling messages are quite small compare with the offloading message body and can be communicated in the dedicated control channel. The proposed distributed algorithm is described in Algorithm \ref{alg:A8}.

\begin{algorithm}[t]
\caption{Distributed Mobile Data Offloading}\label{alg:A8}
\begin{algorithmic}
\STATE Initialize: ${\x}^{0}$,${\y}^{0}$ ${\blambda}^{0}$, $\rho>0$, $\gamma >0$;

\FOR{$k=0,1,\ldots$}

\STATE \COMMENT{Update ${\x}_b$ and ${\y}_a$ for $b = 1,\ldots,B$ and $a = 1,\ldots,A$, \textbf{concurrently}.}

\STATE \COMMENT{\textbf{Base station update, $\forall b$}}

\STATE ${\x}_b^{k+1} \! = \! \argmin_{{\x}_b} \! - \! U_b({\x}_b) \! +\! \frac{\rho}{2}\sum_{a=1}^{A}\Vert x_{ba}\!-\!y_{ab}^{k} - \frac{{\lambda}_{ab}^{k}}{\rho}\Vert_2^{2} \!+\!\frac{1}{2}\Vert{\x}_b \!-\! {\x}_b^{k}\Vert_{\mathbf{P}_i}^{2}$;

\STATE \COMMENT{\textbf{Access point update, $\forall a$}}

\STATE ${\y}_a^{k+1} =  \argmin_{{\y}_b}L_a({\y}_a)+\frac{\rho}{2}\sum_{b=1}^{B}\Vert x_{ba}^{k} - y_{ab} - \frac{{\lambda}_{ab}^{k}}{\rho}\Vert_2^{2} +\frac{1}{2}\Vert{\y}_a - {\y}_a^{k}\Vert_{\mathbf{P}_i}^{2}$;

\STATE \COMMENT{\textbf{SDN controller update}}

\STATE ${\lambda}_{ab}^{k+1} = {\lambda}_{ab}^{k}-\gamma\rho\sum_{b=1}^{B}\sum_{a=1}^{A}(x_{ba}^{k+1}-y_{ab}^{k+1})$;

\ENDFOR
\STATE Output ${\x}$, ${\y}$;
\end{algorithmic}
\end{algorithm}

We consider a wireless access network consists of $B = 5$ base stations and $A = \{5,10\}$ access points coordinated by the SDN controller. The SDN controller will offload mobile data traffic of BSs to APs, and the available capacity of each AP for offloading is $C_a = 10Mbps$. The utility function of BS $b$ is $U_b({\x}_b) = \log({\x}_b^{\top}\mathbf{1} + 1)$, where $\mathbf{1}$ is the all one vector. The cost function of AP $a$ is a linear cost expressed as $L_a({\y}_a) = \theta_a*{\y}_a^{\top}\mathbf{1}$, where $\theta_a > 0$ is the cost coefficient. The value of $\theta_a $ is application specific. During numerical tests, we assume $\theta_a$ is a Gaussian random variable which has a distribution $\mathcal{N}(0,1)$ for simplicity. We perform numerical tests on the offloading decision for one time slot. The simulation result is shown in Figure \ref{fig:convergence}. It shows that the proposed algorithm converges to the optimal objective in a moderate number of iterations when $B=5$ and $A =5$. It takes a longer time for the proposed algorithm to converge when $A = 10$. It indicates that when these are more APs in the access network, it will take a longer time for the SDN controller to coordinate BSs and APs for a consensus on the offloading demand and supply.

\subsection{Other Extensions of ADMM}
\textbf{Decentralized state estimation in smart grid}: Previous work on SCOPF presented in Section. \ref{sec:SCOPF} used direct current (DC) power flow approximation for system state estimation and optimal power flow dispatch. The DC approximation model can provide quick operation instructions for the system. For precise system status monitoring and operation, alternating current (AC) power flow equations are needed
\begin{align}
P_{i} & = \sum_{k=1}^{N}\vert V_i \vert \vert V_k \vert(G_{ik}\cos\theta_{ik} + B_{ik}\sin\theta_{ik}) \\
Q_{i} & = \sum_{k=1}^{N}\vert V_i \vert \vert V_k \vert(G_{ik}\sin\theta_{ik} - B_{ik}\cos\theta_{ik}),
\end{align}
where $P_{i}$ and $Q_{i}$ are real power flow and inactive power flow at bus $i$, respectively. $V_i$ is the voltage magnitude at bus $i$. $G_{ik}$ and $B_{ik}$ are the real and imaginary part of the $(i,k)^{th}$ element of the bus admittance matrix. $\theta_{ik}$ is the voltage phase angle difference between bus $i$ and bus $j$. The problem of state estimation is how to find voltage magnitudes and phase angles given nonlinear equations of real and inactive power flows in the system.

\textbf{Smart meter reading data clustering}: The advanced metering infrastructure (AMI) enables two-way communications with the meter. The smart meters are able to record the consumption of electric energy of each household and send readings to data centers of utility companies for billing and customer service. This provides real time information about electric energy consumption and behaviors of consumers, which can be used for data mining. The smart meters record electric energy consumption of consumers every fifteen minutes, which means that a substantial amount of data are generated daily in the U.S. By investigating those data, we can better understand profiles of consumers to ensure the quality of service, develop targeted electric energy plans, and accurately predict energy consumption of the power system.

\textbf{Efficient air quality monitoring}: The air pollution has been an utmost concern for public health nowadays. In 2012, around seven million people dead worldwide due to the air pollution. However, the existing air-quality monitoring network has very low spatial and temporal coverage, which severely limits its ability to predict air quality and to analyze its impact on environment, climate, and public health. Fortunately, there exists a large amount of diverse data, such as satellite remote sensing data, meteorological data (temperature, wind, pressure, humidity, etc.), and traffic data (volume, speed, congestion) which can be utilized. Instead of solely relying on the traditional monitoring network to provide us the air quality data, many heterogeneous big data sources can be used to develop innovative big data processing methods in air quality research.

\section{Conclusion}
\label{sec:conclusion}
In this paper, we have reviewed several distributed and parallel optimization method based on the ADMM for large scale optimization problems. We have introduced the background of ADMM and described several direct extensions and sophisticated modifications of ADMM from $2$-block to $N$-block settings. We have explained the iterative schemes and convergence properties for each extension/modification. We have illustrated the implantations on large-scale computing facilities, and enumerated several applications of $N$-block ADMM in modern communication networks.

\balance
\bibliographystyle{IEEEtran}
\bibliography{JOC}

\begin{thebibliography}{10}
\providecommand{\url}[1]{#1}
\csname url@samestyle\endcsname
\providecommand{\newblock}{\relax}
\providecommand{\bibinfo}[2]{#2}
\providecommand{\BIBentrySTDinterwordspacing}{\spaceskip=0pt\relax}
\providecommand{\BIBentryALTinterwordstretchfactor}{4}
\providecommand{\BIBentryALTinterwordspacing}{\spaceskip=\fontdimen2\font plus
\BIBentryALTinterwordstretchfactor\fontdimen3\font minus
  \fontdimen4\font\relax}
\providecommand{\BIBforeignlanguage}[2]{{%
\expandafter\ifx\csname l@#1\endcsname\relax
\typeout{** WARNING: IEEEtran.bst: No hyphenation pattern has been}%
\typeout{** loaded for the language `#1'. Using the pattern for}%
\typeout{** the default language instead.}%
\else
\language=\csname l@#1\endcsname
\fi
#2}}
\providecommand{\BIBdecl}{\relax}
\BIBdecl

\bibitem{PS09}
P.~Tseng and S.~Yun, ``A coordinate gradient descent method for non-smooth
  separable minimization,'' \emph{Mathematical Programming}, vol. 117, no.~1,
  pp. 387--423, 2009.

\bibitem{YS09}
Y.~Li and S.~Osher, ``Coordinate descent optimization for l1 minimization with
  applications to compressed sensing: a greedy algorithm,'' UCLA CAM, Tech.
  Rep. Report 09-17, 2009.

\bibitem{Y12}
Y.~Nesterov, ``Efficiency of coordiate descent methods on huge-scale
  optimization problems,'' \emph{SIAM Journal on Optimization}, vol.~22, no.~2,
  pp. 341--362, 2012.

\bibitem{LO08}
L.~Bottou and O.~Bousquet, ``The tradeoffs of large scale learning,'' in
  \emph{Advances in Neural Information Processing Systems}, Vancouver, Canada,
  Dec. 2008.

\bibitem{MMAL10}
M.~Zinkevich, M.~Weimer, A.~Smola, and L.~Li, ``Parallelized stochastic
  gradient descent,'' in \emph{Advances in Neural Information Processing
  Systems}, Vancouver, Canada, Dec. 2010.

\bibitem{FBCS11}
F.~Niu, B.~Recht, C.~Re, and S.~J. Wright, ``Hogwild: A lock-free approach to
  parallelizing stochastic gradient dscent,'' in \emph{Advances in Neural
  Information Processing Systems}, Granada, Spain, Dec. 2011.

\bibitem{CKCSS08}
C.~J. Hsieh, K.~W. Chang, C.~J. Lin, S.~S. Keerthi, and S.~Sundararajan, ``A
  dual coordinate descent method for large-scale linear {SVM},'' in
  \emph{International Conference on Machine Learning}, Helsinki, Finland, Jul.
  2008.

\bibitem{ST13}
S.~Shalev-Shwartz and T.~Zhang, ``Stochastic dual coordinate ascent methods for
  regularized loss minimization,'' \emph{Journal of Machine Learning Research},
  vol.~14, pp. 567--599, 2013.

\bibitem{BT97}
D.~Bertsekas and J.~Tsitsiklis, \emph{Parallel and Distributed Computation:
  Numerical Methods (2nd ed.)}.\hskip 1em plus 0.5em minus 0.4em\relax Belmont,
  MA: Athena Scientific, 1997.

\bibitem{BPCPE10}
S.~Boyd, N.~Parikh, E.~Chu, B.~Peleato, and J.~Eckstein, ``Distributed
  optimization and statistical learning via the alternating direction method of
  multipliers,'' \emph{Foundation and Trends in Machine Learning}, vol.~3,
  no.~1, pp. 1--122, Nov. 2010.

\bibitem{RP14}
R.~M. Freund and P.~Grigas, ``New analysis and results for the frank-wolfe
  method,'' online at http://arxiv.org/abs/1307.0873, 2014.

\bibitem{SMMP13}
S.~Lacoste-Julien, M.~Jaggi, M.~Schmidt, and P.~Pletscher, ``Block-coordinate
  frank-wolfe optimization for structural svms,'' in \emph{International
  Conference on Machine Learning}, Atlanta, GA, Jun. 2013.

\bibitem{RA75}
R.~Glowinski and A.~Marrocco, ``Sur l'approximation par \'el\'ements finis et
  la r\'esolution par p\'enalisation-dualit\'e d'une classe de probl\`emes de
  {D}irichlet non lin\'eaires,'' \emph{Revue Fran\c{c}aise d'Automatique,
  Informatique, Recherche Operationnelle, S\'erie Rouge}, vol. R-2, pp. 41--76,
  1975.

\bibitem{DB76}
D.~Gabay and B.~Mercier, ``A dual algorithm for the solution of nonlinear
  variational problems via finite element approximation,'' \emph{Computers \&
  Mathematics with Applicaions}, vol.~2, no.~1, pp. 17--40, 1976.

\bibitem{YAJWY12}
Y.~Peng, A.~Ganesh, J.~Wright, W.~Xu, and Y.~Ma, ``{RASL}: Robust alignment by
  sparse and lowrank decomposition for linearly correlated images,'' \emph{IEEE
  Transactions on Pattern Analysis and Machine Intelligence}, vol.~34, no.~11,
  pp. 2233--2246, Nov. 2012.

\bibitem{MX11}
M.~Tao and X.~Yuan, ``Recovering low-rank and sparse components of matrices
  from incomplete and noisy observations,'' \emph{SIAM Journal on
  Optimization}, vol.~21, no.~1, pp. 51--87, 2011.

\bibitem{HCB14}
H.~Xu, C.~Feng, and B.~Li, ``Temperature aware workload management in
  geo-distributed datacenters,'' IEEE Transactions on Parallel and Distributed
  Systems, to appear, 2014.

\bibitem{CBYX13}
C.~Chen, B.~He, Y.~Ye, and X.~Yuan, ``The direct extension of {ADMM} for
  multi-block convex minimization problems is not necessarily convergent,''
  preprint, 2013.

\bibitem{MZ12}
M.~Hong and Z.~Luo, ``On the linear convergence of the alternating direction
  method of multipliers,'' online at http://arxiv.org/abs/1208.3922, 2012.

\bibitem{CBYX14}
C.~Chen, B.~He, Y.~Ye, and X.~Yuan, ``The direct extension of {ADMM} for
  multi-block convex minimization problems is not necessarily convergent,''
  Online at http://web.stanford.edu/$\sim$yyye/{ADMM}\_5, 2014.

\bibitem{BXM12}
B.~He, M.~Tao, and X.~Yuan, ``Alternating direction method with {G}aussian back
  substitution for separable convex programming,'' \emph{SIAM Journal of
  Optimization}, vol.~22, no.~2, pp. 313--340, 2012.

\bibitem{MTXMSZ14}
M.~Hong, T.~Chang, X.~Wang, M.~Razaviyayn, S.~Ma, and Z.~Luo, ``A block
  successive upper bound minimization method of multipliers for linearly
  constrained convex optimization,'' online at http://arxiv.org/abs/1401.7079,
  2014.

\bibitem{BLX13}
B.~He, L.~Hou, and X.~Yuan, ``On full jacobian decomposition of the augmented
  lagrangian method for separable convex programming,'' online at
  http://www.optimization-online.org/DB\_HTML/2013/05/3894.html, 2013.

\bibitem{WMZW14}
W.~Deng, M.~Lai, Z.~Peng, and W.~Yin, ``Parallel multi-block {ADMM} with o(1/k)
  convergence,'' online at http://arxiv.org/abs/1312.3040, 2014.

\bibitem{MJM10}
M.~V. Afonso, J.~M. Bioucas-Dias, and M.~A.~T. Figueiredo, ``Fast image
  recovery using variable splitting and constrained optimization,'' \emph{IEEE
  Transactions on Image Processing}, vol.~19, no.~9, pp. 2345--2356, Sep. 2010.

\bibitem{BMX12}
B.~He, M.~Tao, and X.~Yuan, ``Convergence rate and iteration complexity on the
  alternating direction method of multipliers with a substitution procedure for
  separable convex programming,''
  http://www.optimization-online.org/DB\_FILE/2012/09/3611.pdf, 2012.

\bibitem{BMX13}
------, ``On the proximal {J}acobian decomposition of {ALM} for multiple-block
  separable convex minimization problems and its relationship to {ADMM},''
  http://www.optimization-online.org/DB\_FILE/2013/11/4142.pdf, 2013.

\bibitem{HAZ14}
H.~Wang, A.~Banerjee, and Z.~Luo, ``Parallel direction method of multipliers,''
  online at http://arxiv.org/abs/1406.4064, 2014.

\bibitem{NS13}
N.~Parikh and S.~Boyd, ``Proximal algorithms,'' \emph{Foundation and Trends in
  Optimization}, vol.~1, no.~3, pp. 123--231, 2013.

\bibitem{ZMW13}
Z.~Peng, M.~Yan, and W.~Yin, ``Parallel and distributed sparse optimization,''
  in \emph{IEEE Asilomar Conference on Signals, Systems, and Computers},
  Pacific Grove, CA, Nov. 2013.

\bibitem{MMMSI10}
M.~Zaharia, M.~Chowdhury, M.~J. Franklin, S.~Shenker, and I.~Stoica, ``Spark:
  Cluster computing with working sets,'' in \emph{2nd USENIX Conference on Hot
  Topics in Cloud Computing}, Boston, MA, Jun. 2010.

\bibitem{LAWZ13}
L.~Liu, A.~Khodaei, W.~Yin, and Z.~Han, ``A distribute parallel approach for
  big data scale optimal power flow with security constraints,'' in \emph{IEEE
  International Conference on Smart Grid Communications}, Vancouver, Canada,
  Oct. 2013.

\bibitem{SMERS14}
S.~Chakrabarti, M.~Kraning, E.~Chu, R.~Baldick, and S.~Boyd, ``Security
  constrained optimal power flow via proximal message passing,'' in \emph{IEEE
  Power Systems Conference}, Clemson, SC, Mar. 2014.

\bibitem{LXMGZ15}
L.~Liu, X.~Chen, M.~Bennis, G.~Xue, and Z.~Han, ``A distributed admm approach
  for mobile data offloading in software defined network,'' in \emph{IEEE
  Wireless communications and Networking Conference}, New Orleans, LA, March
  2015.

\bibitem{MS14}
M.~Grant and S.~Boyd, ``{CVX}: Matlab software for disciplined convex
  programming, version 2.1,'' \url{http://cvxr.com/cvx}, Mar. 2014.

\end{thebibliography}

\end{document}